\documentclass[12pt,oneside,reqno]{amsart}
\usepackage{mathrsfs}
\usepackage{graphics}
\usepackage{graphicx,color}
\usepackage{amssymb}
\usepackage{subfigure}
\usepackage{enumerate}
\newcommand{\dif}{\mathrm{d}}

\newcommand{\be}{\begin{eqnarray}}
\newcommand{\ee}{\end{eqnarray}}
\newcommand{\ce}{\begin{eqnarray*}}
\newcommand{\de}{\end{eqnarray*}}
\newtheorem{theorem}{Theorem}[section]
\newtheorem{lemma}[theorem]{Lemma}
\newtheorem{remark}[theorem]{Remark}
\newtheorem{definition}[theorem]{Definition}
\newtheorem{proposition}[theorem]{Proposition}
\newtheorem{Example}[theorem]{Example}
\newtheorem{corollary}[theorem]{Corollary}
\def\e{\varepsilon}

\def\g{\gamma}

\def\[{{\Big[}}
\def\]{{\Big]}}
\def\<{{\langle}}
\def\>{{\rangle}}
\def\({{\Big(}}
\def\){{\Big)}}

\def\no{\nonumber}
\def\bt{\begin{theorem}}
\def\et{\end{theorem}}
\def\bl{\begin{lemma}}
\def\el{\end{lemma}}
\def\br{\begin{remark}}
\def\er{\end{remark}}
\def\bx{\begin{Example}}
\def\ex{\end{Example}}
\def\bd{\begin{definition}}
\def\ed{\end{definition}}
\def\bp{\begin{proposition}}
\def\ep{\end{proposition}}
\def\bc{\begin{corollary}}
\def\ec{\end{corollary}}

\def\mE{{\mathbb E}}

\def\mG{{\mathbb G}}
\def\mH{{\mathbb H}}

\def\mN{{\mathbb N}}

\def\mP{{\mathbb P}}

\def\mR{{\mathbb R}}

\def\mU{{\mathbb U}}
\def\mV{{\mathbb V}}

\def\geq{\geqslant}
\def\leq{\leqslant}

\begin{document}

\allowdisplaybreaks

\title{Support theorems for degenerate stochastic
differential equations with jumps and applications*}

\author{Huijie Qiao$^1$ and Jiang-Lun Wu$^2$}



\thanks{*This work was partly supported by NSF of China (No. 11001051, 11371352, 11671083) and China Scholarship Council under Grant No. 201906095034.}


\date{}

\dedicatory{1. School of Mathematics, Southeast University\\
Nanjing, Jiangsu 211189,  China\\
hjqiaogean@seu.edu.cn\\
2. Department of Mathematics, Computational Foundry, Swansea University\\
Bay Campus, Swansea SA1 8EN, UK\\
j.l.wu@swansea.ac.uk}

\begin{abstract}
In the paper, we are concerned with degenerate stochastic differential equations with jumps. Firstly, we establish two support theorems for the solutions of the degenerate stochastic differential equations,
under different (sufficient) conditions. Secondly, we apply one of our support theorems to a class of degenerate stochastic evolution equations (i.e., infinite-dimensional stochastic differential equations)
with jumps to get a characterisation of path-independence for the densities of their Girsanov transformations.
\end{abstract}

\maketitle \rm

\section{Introduction}

Support theorems for stochastic differential equations (SDEs) are referred to as the supports for the laws (or the distributions) of their solutions, or equivalently, are that the state spaces for the solutions are characterized under certain topologies. Support theorems for diffusion processes governed by SDEs were initiated in the two seminar papers Stroock-Varadhan \cite{sv,sv1}. Since then, there have been a lot of results on this topic. Let us recall some works related to ours here. An extension of the celebrated Stroock-Varadhan support theorem to SDEs with jumps is due to Simon in \cite{s} who established support theorems for a class of SDEs driven by Poisson random measures. Later, Fournier \cite{nf} investigated a class of parabolic stochastic partial differential equations driven by space-time Gaussian white noises and independent Poisson measures, and proved that the supports of their solution distributions are characterized as the closures of sets of weak solutions for the corresponding (deterministic) partial differential equations of parabolic type. However, for a fixed time, the supports for the distributions of solutions of SDEs with jumps are less considered. 

Besides, SDEs with degenerate coefficients and jumps attract much more attentions in recent years, see e.g. \cite{zxc0,zxc} and references therein. Thus, to study the support property for degenerate SDEs with jumps is interesting and potentially useful. In the present paper, we are concerned with the supports of a class of degenerate SDEs driven by Brownian motions and independent compensated Poisson random measures for a fixed time. We aim to establish support theorems for these SDEs. We then utilise our result to derive a characterisation theorem for path-independent property of Girsanov transformation for degenerate stochastic evolution equations with jumps.
More precisely, in the paper, we will prove a theorem on supports of distributions for SDEs under certain general assumptions. And then by strengthening the conditions such that SDEs have densities, and applying the fact that existence of the density for an SDE implies full support, we obtain the other support theorem. Next, we make use of one support theorem to a problem on degenerate stochastic evolution equations with jumps on Hilbert spaces. By some deduction, path-independence for the densities of their Girsanov transformations is characterized.

It is worthwhile to mentioning our previous results. In \cite{qw1, qw2}, we showed that these densities of Girsanov transformations for SDEs with jumps and stochastic evolution equations with jumps are path-independent, when the coefficients of their continuous diffusion terms are non-degenerate. And then the second named author and B. Wu \cite{ww} only mentioned that these densities of Girsanov transformations for degenerate stochastic differential equations are path-independent. Here in this paper, we permit that their continuous diffusion coefficients are degenerate, and give some concrete conditions and detailed proof.

This rest of the paper is organized as follows. In Section \ref{suth}, we prove two support theorems for SDEs with jumps under different (sufficient) conditions. Section \ref{app} is devoted to applying a support theorem to a problem on degenerate stochastic evolution equations with jumps. Finally, we obtain path-independence for the density of the  associated Girsanov transformation.

The following convention will be used throughout the paper: $C$ with or without indices will denote different positive constants (depending on the indices) whose values may change from one place to another.

\section{Support theorems}\label{suth}

In the section, we will prove two support theorems for stochastic differential equations with jumps applied in the next section.

Let $(\Omega,\mathscr{F},\mP;(\mathscr{F}_t)_{t\geq 0})$ be a complete filtered probability space. Let $(\mU, \mathscr{U}, \|\cdot\|_{\mU})$ be a finite dimensional normed space. Let $\nu$ be a $\sigma$-finite measure defined on $(\mU,\mathscr{U})$. We fix $\mU_0\in\mathscr{U}, \mU_0\subset\mU-\{0\}$ with $\nu(\mU\setminus\mU_0)<\infty$ and $\int_{\mU_0}\|u\|_{\mU}^2\,\nu(\dif u)<\infty$. And then we construct an integer-valued $(\mathscr{F}_t)_{t\geq 0}$-Poisson random measure $N(\dif t, \dif u)$ on $(\Omega,\mathscr{F},\mP;(\mathscr{F}_t)_{t\geq 0})$ with the intensity $\dif t\nu(\dif u)$. Set
$$\tilde{N}(\dif t,\dif u):=N(\dif t,\dif u)-\dif t\nu(\dif u)$$
and then $\tilde{N}(\dif t,\dif u)$ is the compensated $(\mathscr{F}_t)_{t\geq 0}$-predictable martingale
measure of $N(\dif t, \dif u)$. And then let $\{B_t\}$ be a $d$-dimensional $(\mathscr{F}_t)_{t\geq 0}$-Brownian motion, which is independent of  $N(\dif t, \dif u)$. Fix $T>0$ and consider the following SDE with jumps on $\mR^d$
\be\left\{\begin{array}{l}
\dif Z_t=\xi(Z_t)\dif t+\eta(Z_t)\dif B_t+\int_{\mU_0}\zeta(Z_{t-},u)\tilde{N}(\dif t, \dif u), \qquad t\in(0,T],\\
Z_0=\g.
\end{array}
\right.
\label{Eq1}
\ee
The coefficients $\xi: \mR^d\mapsto\mR^d$, $\eta: \mR^d\mapsto\mR^{d\times d}$ and
$\zeta: \mR^d\times\mU_0\mapsto\mR^d$ are all Borel measurable. And $\g$ is a $\mathscr{F}_0$-measurable random variable with $\mE|\g|^2<\infty$.

\begin{enumerate}[\bf{Assumption 1.}]
\item
\end{enumerate}
\begin{enumerate}[{\bf (i)}]
\item There exists a constant $L_1>0$ such that for any $z, z_1, z_2\in\mR^d$ and $u, u_1, u_2\in\mU_0$,
$$
|\xi(z_1)-\xi(z_2)|+\|\eta(z_1)-\eta(z_2)\|\leq L_1|z_1-z_2|,
$$
where $\|\cdot\|$ stands for the Hilbert-Schmidt norm of a matrix, and
\ce
&&|\zeta(z_1,u)-\zeta(z_2,u)|\leq L_1|z_1-z_2|\|u\|_{\mU},\\
&&|\zeta(z,u_1)-\zeta(z,u_2)|\leq L_1(1+|z|)\|u_1-u_2\|_{\mU}.
\de
\item There exists a constant $L_2>0$ such that for any $z\in\mR^d$
\ce
|\xi(z)|^2+\|\eta(z)\|^2+\int_{\mU_0}|\zeta(z,u)|^2\nu(\dif u)
\leq L_2(1+|z|^2).
\de
\end{enumerate}

\vspace{3mm}

Under {\bf Assumption 1}, it is known that there exists a unique strong solution $\{Z_t\}_{t\in[0,T]}$ to Eq.(\ref{Eq1}) which is a Markov process with c\`adl\`ag paths, see, e.g., \cite[Theorem 6.2.3 and Theorem 6.4.5]{da}. Next, we define the support for a random variable and then study the support of $Z_t$ for $t\in[0,T]$.

\bd\label{supset}
Let $\mV$ be a metric space with the metric $\rho$. The support of a $\mV$-valued random variable $v$ is defined to be
$$
{\rm supp}(v):=\{z\in\mV | (\mP\circ v^{-1})(B(z,r))>0, ~\mbox{for all}~ r>0\},
$$
where $B(z,r):=\{y\in\mV | \rho(z,y)<r\}$.
\ed

\vspace{3mm}

\begin{enumerate}[\bf{Assumption 2.}]
\item For any $z\in\mR^d$ and any open ball $B\subset\mR^d$, there exists a point $u\in {\rm supp}(\nu^{\mU_0})$, where $\nu^{\mU_0}$ is the restriction of $\nu$ to $\mU_0$,  such that $\zeta(z,u)\in B$.
\end{enumerate}

\vspace{3mm}

Now, we state and prove the first main result of this section.
\bt\label{supth}
Suppose that ${\rm supp}(\g)=\mR^d$ and $\xi, \eta, \zeta$ satisfy {\bf Assumption 1} and {\bf Assumption 2}. Then ${\rm supp}(Z_t)=\mR^d$ for $t\in[0,T]$.
\et
\begin{proof}
By Definition \ref{supset}, it is sufficient to prove that for any $t\in[0, T], a\in\mR^d$ and $r>0$
$$
\mP(Z_t\in B(a, r))>0.
$$
And then we fix $t, a, r$ and prove the above inequality with the help of two auxiliary processes.

{\bf Step 1.} For any subset $U\subset\mU_0$ with $U\in\mathscr{U}, \nu(U)<\infty$ and $L_1\|u\|_{\mU}<1$ for $u\in U$, we introduce the first auxiliary equation
\be
Z_t^U=\g+\int_0^t\xi(Z_s^U)\dif s+\int_0^t\int_{U}\zeta(Z^U_{s-},u)\tilde{N}(\dif s, \dif u).
\label{auxpr}
\ee
Note that there is no continuous diffusion term in the above equation. This is because we permit that Eq.(\ref{Eq1}) is degenerate and then the continuous diffusion term is not needed. Under {\bf Assumption 1}, it follows from \cite[Theorem 6.2.3]{da} that Eq.(\ref{auxpr}) has a unique solution denoted as $Z^U$. Thus, Lemma \ref{firest} below admits us to obtain that for a.s. $\omega\in\Omega$ and any $\e\in(0,r)$, there exists a $n\in\mN$ such that
\be
\sup\limits_{0\leq s\leq t_n}|Z_s-Z_s^U|<\e/2.
\label{auxfirest1}
\ee
Fix $t_n$ and $\e$.

{\bf Step 2.} Let $\{s_i\}$ be a positive sequence such that $s_i\uparrow\infty$, and $\{u_i\}$ be a sequence in the support of $\nu^U$,  the restriction of $\nu$ to $U$. And let $\mG^U$ be the collection of the above sequences pair $\{s_i\}, \{u_i\}$. For any $g\in\mG^U$, we introduce the second auxiliary equation
\be
Z^{g, U}_t=\g+\int_0^t\left[\xi(Z^{g, U}_s)-\int_{U}\zeta(Z^{g, U}_s,u)\nu(\dif u)\right]\dif s+\sum_{i:s_i\leq t}\zeta(Z^{g, U}_{s_i-},u_i).
\label{auxeq2}
\ee
Under {\bf Assumption 1}, it holds that Eq.(\ref{auxeq2}) has a unique solution denoted as $Z^{g, U}$. So, by {\bf Assumption 2}, we know that for the open ball $B(a, r-\e)$, there exist $s_1>0, s_i>t_n, i=2, 3, \dots$ and $u_1\in {\rm supp}(\nu^U)$ such that
\be
|Z^{g, U}_{t_n}-a|<r-\e.
\label{auxfirest2}
\ee
This may be possible if $s_1$ is taken enough small so that $s_1\leq t_n$ and $Z^{g, U}_{t_n}\in B(a, r-\e)$. Fix $\{s_i\}$ and $u_1$.

{\bf Step 3.} We study the relationship between $Z^U$ and $Z^{g, U}$. Set
\ce
\chi_t:=\int_0^t\int_{U}u\tilde{N}(\dif s, \dif u), \quad \Delta\chi_t:=\chi_t-\chi_{t-},
\quad D:=\{t\in[0,\infty), \Delta\chi_t\in U\},
\de
and then it follows from $\nu(U)<\infty$ that $D$ is a discrete set in $[0,\infty)$ a.s.. Let $0<\tau_1<\tau_2<\cdots<\tau_n<\cdots$ be the enumeration of all elements in $D$. Besides, we  take any $u_2\in supp(\nu^U)$. For any $\e'>0$, set
\ce
&&A_1:=\{\omega\in\Omega: 0<s_1-\tau_1<\e', \|u_1-\Delta\chi_{\tau_1}\|_{\mU}<\e'\},\\
&&A_2:=\{\omega\in\Omega: 0<s_2-s_1-(\tau_2-\tau_1)<\e', \|u_2-\Delta\chi_{\tau_2}\|_{\mU}<\e'\},
\de
and then it follows from independence for increments of $N(\dif t, \dif u)$ that $\mP(A_1\cap A_2)>0$. Thus, by Lemma \ref{secest} below it holds that for the above $\e$, there exists an $\e'>0$ such that
\be
\sup\limits_{0\leq s\leq t_n}|Z_s^U-Z^{g, U}_s|<\e/2,
\label{auxfirest3}
\ee
on $A_1\cap A_2$.

{\bf Step 4.} Combining (\ref{auxfirest1}) (\ref{auxfirest2}) with (\ref{auxfirest3}), we obtain that
\ce
|Z_s-a|\leq|Z_s-Z_s^U|+|Z_s^U-Z^{g, U}_s|+|Z^{g, U}_s-a|<\e/2+\e/2+r-\e=r, \quad s\in(0, t_n],
\de
on $A_1\cap A_2$. Thus $\mP(|Z_s-a|<r)>0$ for $s\in(0, t_n]$. If $t\leq t_n$, the proof is over; if $t>t_n$, by the Markov property, we still can obtain $\mP(|Z_t-a|<r)>0$. The proof is completed.
\end{proof}

\bl\label{firest}
Under {\bf Assumption 1}, there exists a positive (nonrandom) sequence $\{t_n\}$ decreasing to $0$ such that
\ce
\lim_{n\rightarrow\infty}\sup_{0\leq s\leq t_n}|Z_s-Z_s^U|=0, \quad a.s. \,\mP.
\de
\el
\begin{proof}
Firstly, we compute $Z_t-Z_t^U$ for $t\in[0,T]$. By (\ref{Eq1}) and (\ref{auxpr}), it holds that
\ce
Z_t-Z_t^U&=&\int_0^t\left(\xi(Z_s)-\xi(Z_s^U)\right)\dif s+\int_0^t\eta(Z_s)\dif B_s\\
&&+\int_0^t\int_{U}\left(\zeta(Z_{s-},u)-\zeta(Z^U_{s-},u)\right)\tilde{N}(\dif s, \dif u)\\
&&+\int_0^t\int_{\mU_0\setminus U}\zeta(Z_{s-},u)\tilde{N}(\dif s, \dif u).
\de
And by the Burkholder-Davis-Gundy inequality and the H\"older inequality, one can have that
\ce
\mE\left(\sup_{0\leq s\leq t}|Z_s-Z_s^U|^2\right)&\leq&4t\mE\int_0^t\left|\xi(Z_s)-\xi(Z_s^U)\right|^2\dif s+16\mE\int_0^t|\eta(Z_s)|^2\dif s\\
&&+16\mE\int_0^t\int_{U}\left|\zeta(Z_{s-},u)-\zeta(Z^U_{s-},u)\right|^2\nu(\dif u)\dif s\\
&&+16\mE\int_0^t\int_{\mU_0\setminus U}|\zeta(Z_{s-},u)|^2\nu(\dif u)\dif s.
\de
Moreover, by {\bf Assumption 1}, we obtain that
\be
\mE\left(\sup_{0\leq s\leq t}|Z_s-Z_s^U|^2\right)&\leq&4L^2_1(t+4\nu(U))\int_0^t\mE\left(\sup_{0\leq s\leq r}|Z_s-Z_s^U|^2\right)\dif r\no\\
&&+16\mE\int_0^tL_2(1+|Z_s|^2)\dif s.
\label{mides}
\ee

To estimate the last term in (\ref{mides}), we observe Eq.(\ref{Eq1}). By similar deduction to above, one can get that
\ce
\mE\left(\sup_{0\leq s\leq t}|Z_s|^2\right)&\leq&4\mE|\g|^2+4t\mE\int_0^t\left|\xi(Z_s)\right|^2\dif s+16\mE\int_0^t|\eta(Z_s)|^2\dif s\\
&&+16\mE\int_0^t\int_{\mU_0}|\zeta(Z_{s-},u)|^2\nu(\dif u)\dif s,
\de
and furthermore by {\bf Assumption 1}
\ce
\mE\left(\sup_{0\leq s\leq t}|Z_s|^2\right)&\leq&4\mE|\g|^2+4(t+4)L_2\int_0^t\mE\left(1+|Z_s|^2\right)\dif s\\
&\leq&4\mE|\g|^2+4(t+4)tL_2+4(t+4)L_2\int_0^t\mE\left(\sup_{0\leq s\leq r}|Z_s|^2\right)\dif r.
\de
Thus, the Gronwall inequality admits us to have that
\be
\mE\left(\sup_{0\leq s\leq t}|Z_s|^2\right)\leq C,
\label{mides2}
\ee
where the constant $C>0$ depends on $\mE|\g|^2, T, L_2$.

Next, combining (\ref{mides}) with (\ref{mides2}), we get that
\ce
\mE\left(\sup_{0\leq s\leq t}|Z_s-Z_s^U|^2\right)&\leq&16L_2(C+1)t+4L^2_1(T+4\nu(U))\int_0^t\mE\left(\sup_{0\leq s\leq r}|Z_s-Z_s^U|^2\right)\dif r.
\de
Based on the Gronwall inequality, it holds that
\ce
\mE\left(\sup_{0\leq s\leq t}|Z_s-Z_s^U|^2\right)\leq C(e^{Ct}-1),
\de
where the constant $C>0$ depends on $\mE|\g|^2, T, L_1, L_2$. Set $t_n:=C^{-1}\ln(1+2^{-n})$, for $n=1, 2, \cdots$, and then
\ce
\mE\left(\sup_{0\leq s\leq t_n}|Z_s-Z_s^U|^2\right)\leq C2^{-n}.
\de
Thus, there exists a subsequence still denoted as $\{t_n\}$ such that
\ce
\lim_{n\rightarrow\infty}\sup_{0\leq s\leq t_n}|Z_s-Z_s^U|=0, \quad a.s. \mP.
\de
The proof is completed.
\end{proof}

\bl\label{secest}
Under {\bf Assumption 1}, for the above $\e$, there eixsts an $\e'>0$ such that
$$
\sup_{0\leq s\leq t_n}|Z_s^U-Z^{g, U}_s|<\e/2,
$$
on $A_1\cap A_2$.
\el
\begin{proof}
By (\ref{auxpr}) and (\ref{auxeq2}), it holds that for $0\leq t\leq t_n$
\ce
Z_t^U-Z^{g, U}_t&=&\int_0^t\left[\xi(Z^U_s)-\xi(Z^{g, U}_s)\right]\dif s-\int_0^t\int_{U}\left[\zeta(Z^{U}_s,u)-\zeta(Z^{g, U}_s,u)\right]\nu(\dif u)\dif s\\
&&+\zeta(Z^{U}_{\tau_1-},\Delta\chi_{\tau_1})-\zeta(Z^{g, U}_{s_1-},u_1),
\de
and
\ce
|Z_t^U-Z^{g, U}_t|&\leq&\int_0^t\left|\xi(Z^U_s)-\xi(Z^{g, U}_s)\right|\dif s+\int_0^t\int_{U}\left|\zeta(Z^{U}_s,u)-\zeta(Z^{g, U}_s,u)\right|\nu(\dif u)\dif s\\
&&+|\zeta(Z^{U}_{\tau_1-},\Delta\chi_{\tau_1})-\zeta(Z^{U}_{s_1-},\Delta\chi_{\tau_1})|+|\zeta(Z^{U}_{s_1-},\Delta\chi_{\tau_1})-\zeta(Z^{U}_{s_1-},u_1)|\\
&&+|\zeta(Z^{U}_{s_1-},u_1)-\zeta(Z^{g, U}_{s_1-},u_1)|.
\de
So, {\bf Assumption 1} admits us to obtain that
 \ce
\sup_{0\leq s\leq t_n}|Z_s^U-Z^{g, U}_s|&\leq&L_1\int_0^{t_n}\sup_{0\leq s\leq r}|Z_s^U-Z^{g, U}_s|\dif r\\
&&+L_1\left(\int_{U}\|u\|_{\mU}\nu(\dif u)\right)\int_0^{t_n}\sup_{0\leq s\leq r}|Z_s^U-Z^{g, U}_s|\dif r\\
&&+L_1\|\Delta\chi_{\tau_1}\|_{\mU}|Z^{U}_{\tau_1-}-Z^{U}_{s_1-}|+L_1(1+|Z^{U}_{s_1-}|)\|\Delta\chi_{\tau_1}-u_1\|_{\mU}\\
&&+L_1\|u_1\|_{\mU}|Z^{U}_{s_1-}-Z^{g, U}_{s_1-}|\\
&\leq&L_1\left(1+\int_{U}\|u\|_{\mU}\nu(\dif u)\right)\int_0^{t_n}\sup_{0\leq s\leq r}|Z_s^U-Z^{g, U}_s|\dif r\\
&&+L_1\|\Delta\chi_{\tau_1}\|_{\mU}|Z^{U}_{\tau_1-}-Z^{U}_{s_1-}|+L_1(1+|Z^{U}_{s_1-}|)\|\Delta\chi_{\tau_1}-u_1\|_{\mU}\\
&&+L_1\|u_1\|_{\mU}\sup_{0\leq s\leq t_n}|Z_s^U-Z^{g, U}_s|.
\de
By {\bf Assumption 1}, the definition of $A_1, A_2$ and the Gronwall inequality, we know that there exists a constant $C>0$ such that
$$
|Z^{U}_{s_1-}|<C, \quad |Z^{U}_{\tau_1-}-Z^{U}_{s_1-}|<C\e'.
$$
Thus, on $A_1$ one have
\ce
\sup_{0\leq s\leq t_n}|Z_s^U-Z^{g, U}_s|&\leq&\frac{L_1\left(1+\int_{U}\|u\|_{\mU}\nu(\dif u)\right)}{1-L_1\|u_1\|_{\mU}}\int_0^{t_n}\sup_{0\leq s\leq r}|Z_s^U-Z^{g, U}_s|\dif r\\
&&+\frac{L_1(\e'+\|u_1\|_{\mU})C}{1-L_1\|u_1\|_{\mU}}\e'+\frac{L_1(1+C)}{1-L_1\|u_1\|_{\mU}}\e'.
\de
The Gronwall inequality admits us to obtain that
$$
\sup_{0\leq s\leq t_n}|Z_s^U-Z^{g, U}_s|\leq C\e'.
$$
Taking $C\e'=\e/2$, one can attain that
$$
\sup_{0\leq s\leq t_n}|Z_s^U-Z^{g, U}_s|<\e/2,
$$
on $A_1\cap A_2$. The proof is completed.
\end{proof}

\vspace{6mm}

Next, we want to strengthen those conditions in {\bf Assumptions 1 and 2} and give the other support theorem.

\vspace{3mm}

\begin{enumerate}[\bf{Assumption 3.}]
\item
\end{enumerate}
\begin{enumerate}[{\bf (i)}]
\item The distribution of $\g$ is absolutely continuous with respect to the Lebesgue measure on $\mR^d$.
\item $\xi$ and $\eta$ are $3$-times differentiable with bounded derivatives of all order between $1$ and $3$.
\item For any $u\in\mU_0$, $\zeta(\cdot,u)$ is $3$-times differentiable, and
\ce
\zeta(0,\cdot)&\in&\bigcap_{2\leq q<\infty}L^q(\mU_0,\nu)\\
\sup\limits_{x}|\partial^r_x\zeta(x,\cdot)|&\in&\bigcap_{2\leq q<\infty}L^q(\mU_0,\nu),
\quad 1\leq r\leq3,
\de
where $\partial^r_x\zeta(x,\cdot)$ stands for $r$ order partial derivative of $\zeta(x,\cdot)$
with respect to $x$.
\end{enumerate}

\vspace{3mm}

Under {\bf Assumption 3}, by \cite[Theorem 2-14, p.11]{bgj}, it holds that Eq.(\ref{Eq1}) has a unique solution
which is still denoted by $Z_t$ and the distribution of $Z_t$ possesses a density. Our second main result of this section states
as follows.

\bt\label{sesuth}
Suppose that {\bf Assumption 3} is satisfied. Then ${\rm supp}(Z_t)=\mR^d$ for $t\in[0,T]$.
\et
\begin{proof}
By \cite[Theorem 2-14, p.11]{bgj}, we know that the distribution of $Z_t$ is absolutely continuous with respect to the Lebesgue measure on $\mR^d$ and clearly, the support of the Lebesgue measure is the whole $\mR^d$, the support for the distribution of $Z_t$ is also $\mR^d$. The proof is completed.
\end{proof}

\section{Application}\label{app}

In the section, we will apply Theorem \ref{supth} to a problem on stochastic evolution equations with jumps on (separable) Hilbert spaces.

Let us begin with some notions and notations. Let $\mH$ be a Hilbert space with the inner product $\langle\cdot,\cdot\rangle_{\mathbb{H}}$ and the norm $\|\cdot\|_{\mathbb{H}}$. Let $L(\mH)$ be the set of
all bounded linear operators $L:$ $\mH\rightarrow\mH$ and $L_{HS}(\mH)$ be the collection of all Hilbert-Schmidt
operator $L:\mH\rightarrow\mH$ equipped with the Hilbert-Schmidt norm $\|\cdot\|_{HS}$.

Let $A$ be a linear, unbounded, negative definite and self-adjoint operator on $\mH$ and $D(A)$ be the domain of the operator $A$. Let $\{e^{tA}\}_{t\ge0}$ be the contraction $C_0$-semigroup generated by $A$. Let $L_{A}(\mH)$ be the collection of all densely defined closed linear operators $(L,D(L))$ on $\mH$ so that $e^{t A}L$ can
extend uniquely to a Hilbert-Schmidt operator still denoted by $e^{t A}L$ for any $t>0$. And then $L_{A}(\mH)$, endowed with the $\sigma$-algebra induced by $\{L\to\langle e^{tA}L x,y\rangle_{\mH}\mid t>0,x,y\in \mH\}$, becomes a measurable space.

 Let
$\{\beta^i, i\in\mN\}$ be a family of mutually independent one-dimensional Brownian motions on
$(\Omega,\mathscr{F},\mP;(\mathscr{F}_t)_{t\geq 0})$.  So, we can construct a cylindrical Brownian motion on $\mathbb{H}$ by
\begin{equation*}
W_t:=\sum_{i=1}^\infty\beta^i_te_i, \quad \,\, t\in[0,\infty),
\end{equation*}
where $\{e_i, i\in\mN\}$ is a complete orthonormal basis for $\mathbb{H}$ which will be specified later.  It can be justified
that the covariance operator of the cylindrical Brownian motion $W$ is the identity operator $I$ on $\mathbb{H}$. It is worthwhile to mention that $W$ is not a process on $\mathbb{H}$. However, $W$ can be realized as a continuous process on an enlarged Hilbert space $\tilde{\mathbb{H}}$, the completion of $\mathbb{H}$ under the inner product
$$
\<x,y\>_{\tilde{\mathbb{H}}}:=\sum\limits_{i=1}^\infty 2^{-i}\<x,e_i\>_{\mH}\<y,e_i\>_{\mH}, \quad x, y\in\mathbb{H}.
$$

Next, we introduce a type of jump measures. Let $\lambda:\mathbb{U}\to(0,1)$ be a measurable
function. Then, by Theorem I.8.1 of \cite{iw}, we can construct an integer-valued random measure on $[0,\infty)\times\mathbb{U}$
$$N_{\lambda}:\mathscr{B}([0,\infty))\times\mathscr{U}\times\Omega\to\mathbb{N}_0:=
\mathbb{N}\cup\{0\}$$
with the predictable compensator $\lambda(u)\dif t\nu(\dif u)$:
$$\mathbb{E} N_{\lambda}(\dif t, \dif u,\cdot)=\lambda(u)\dif t\nu(\dif u).$$
Set
$$\tilde{N}_\lambda(\dif t,\dif u):=N_\lambda(\dif t,\dif u)-\lambda(u)\dif t\nu(\dif u),$$
and then $\tilde{N}_\lambda(\dif t,\dif u)$  is the associated compensated martingale measure of $N_{\lambda}(\dif t, \dif u)$. Moreover,
we assume that $W_t, N_{\lambda}(\dif t,\dif u)$ are mutually independent.

Now consider the following stochastic evolution equation with jumps on $\mH$
\begin{equation}
\left\{ \begin{aligned}
         d X_t&=\{A X_t+b(X_t)\}dt+\sigma(X_t)dW_t+\int_{\mU_0}f(X_{t-},u)\tilde{N_\lambda}(\dif t, \dif u),\ \ \ 0<t\leq T,\\
                  X_0&=\Gamma,
                          \end{aligned} \label{3} \right.
                          \end{equation}
where $b:\mH\rightarrow\tilde{\mH}$, $\sigma:\mH\rightarrow L_{A}(\mH)$ and $f:\mH\times\mU_0\rightarrow\tilde{\mH}$ are all Borel measurable mappings, 
and $\Gamma$ is a $\mathscr{F}_0$-measurable random variable with $\mE|\Gamma|^2<\infty$ and ${\rm supp}(\Gamma)=\mH$. Set $\|x\|_{\mH}=\infty, x\notin\mH$. For $b, \sigma, f$, we make the following assumption.

\vspace{3mm}

\begin{enumerate}[\bf{Assumption 4.}]
\item
\end{enumerate}
\begin{enumerate}[{\bf (i)}]
\item There exists an integrable function $L_b:(0,T]\to(0,\infty)$ such that
$$
\|e^{s A}(b(x)-b(y))\|^2_{\mH}\leq L_b(s)\|x-y\|^2_{\mH}, \ \ \ s\in(0,T], x,y\in\mH,
$$
and
$$
\int^T_0\|e^{s A}b(0)\|^2_{\mH}\dif s<\infty.
$$
\item There exists an integrable function $L_{\sigma}:(0,T]\to(0,\infty)$ such that for $\forall s\in(0,T]$ and $\forall x,y\in \mH$
$$
\|e^{s A}\left(\sigma(x)-\sigma(y)\right)\|^2_{HS}\leq L_{\sigma}(s)\|x-y\|^2_{\mH}$$
and
$$
\int^T_0\|e^{s A}\sigma(0)\|^2_{HS}\dif s<\infty.
$$
\item There exists an integrable function $L_f:[0,T]\to(0,\infty)$ such that
\ce
&&\|e^{s A}(f(x,u)-f(y,u))\|^2_{\mH}\leq L_f(s)\|u\|^2_{\mU}\|x-y\|^2_{\mH}, \ \ \ s\in[0,T], u\in\mU_0, x,y\in \mH,\\
&&\|e^{s A}(f(x,u_1)-f(x,u_2))\|^2_{\mH}\leq L_f(s)(1+\|x\|_{\mH})^2\|u_1-u_2\|^2_{\mU}, \ \ \ u_1, u_2\in\mU_0,
\de
and
$$
\int_{\mU_0}\|e^{s A}f(x,u)\|^2_{\mH}\lambda(u)\nu(\dif u)\leq L_f(s)(1+\|x\|_{\mH})^2.
$$
\end{enumerate}

\vspace{3mm}

Under {\bf Assumption 4}, \cite[Theorem 3.2]{qw2} admits us to obtain that Eq(\ref{3}) has a unique mild solution, denoted by $X_t$.

\vspace{3mm}

\begin{enumerate}[\bf{Assumption 5.}]
\item The operator $-A$ has the following eigenvalues
$$0<\lambda_1<\lambda_2<\dots <\lambda_j<\ldots$$
counting multiplicities.
\end{enumerate}

\vspace{3mm}

The complete orthonormal basis $\{e_j\}_{j\in\mathbb{N}}$ of $\mH$ is taken as the eigen-basis of $-A$ throughout the rest of the paper. Let ${\mH}_n$ be the space spanned by $e_1,\cdots,e_n$ for $n\in\mN$. Define
\ce
\pi_n:\mH\rightarrow {\mH}_n,\qquad \pi_nx:=\sum_{i=1}^n\langle x,e_i\rangle_{\mH}e_i, \quad x\in\mH,
\de
and then $\pi_n$ is the orthogonal project operator from $\mH$ to ${\mH}_n$.

\vspace{3mm}

\begin{enumerate}[\bf{Assumption 6.}]
\item For any $n\in\mN, z\in\mH_n$ and any open ball $B\subset\mH_n$, there exists a point $u\in {\rm supp}(\nu^{\mU_0})$ such that $\pi_nf(z,u)\in B$.
\end{enumerate}

\vspace{3mm}

In the following, we give out a support theorem under these assumptions.

\bl\label{supthe}
Under {\bf Assumptions 4-6}, it holds that ${\rm supp}(X_t)=\mH$ for $t\in[0,T]$.
\el
\begin{proof}
Since ${\rm supp}(X_t)\subset\mH$, it is sufficient only to show that ${\rm supp}(X_t)\supset\mH$. Furthermore, by Definition \ref{supset},
we only need to prove that for any $h\in\mH$ and $r>0$,
$$
\mP\{\|X_t-h\|_{\mH}<r\}>0,
$$
or equivalently,
$$
\mP\{\|X_t-h\|_{\mH}\geq r\}<1.
$$

On one hand,  put $A_n:=A\mid_{{\mH}_n}, b_n:=\pi_n b$, $\sigma_n:=\pi_n\sigma$ and $f_n:=\pi_n f$. Thus, we obtain the following SDE with jumps in ${\mH}_n$
\begin{equation}\label{8}
\begin{cases}
d X^n_t=\{A_nX^n_t+b_n(X^n_t)\}dt+\sigma_n(X^n_t)dW_t+\int_{\mU_0}f_n(X^n_{t-},u)\tilde{N}_\lambda(\dif t, \dif u),\\
X^n(0)=\pi_n\Gamma.
\end{cases}
\end{equation}
It is easy to see that Eq.(\ref{8}) is similar to Eq.(\ref{Eq1}). Moreover, it follows from {\bf Assumptions 4-6} that the coefficients $b_n$, $\sigma_n$ and $f_n$ satisfy {\bf Assumptions 1-2}. Thus, by Theorem \ref{supth}, it holds that Eq.(\ref{8}) has a unique solution which is denoted by $X_t^n$ with ${\rm supp}(X_t^n)={\mH}_n$. So, Definition \ref{supset} admits us to obtain that for any small $0<\e<r$ and $0<\eta<1$
$$
\mP\{\|X^n_t-\pi_n h\|_{\mH}\geq r-\e\}<1-\eta.
$$

On the other hand, it follows from \cite[Lemma 3.3]{qw2} that
$$
\lim_{n\rightarrow\infty}\mathbb{E}\|X_t^n-X_t\|_{\mH}^2=0, \quad t\in[0,T].
$$
Thus, by the Chebyshev inequality we know that
there exists a $N\in\mN$ such that for $n>N$,
$$
\mP\{\|X_t-X^n_t\|_{\mH}\geq \e/2\}<\eta/2, \quad \mP\{\|\pi_nh-h\|_{\mH}\geq \e/2\}<\eta/2.
$$

Finally, based on these inequalities, it holds that
\ce
\mP\{\|X_t-h\|_{\mH}\geq r\}&\leq&\mP\{\|X_t-X^n_t\|_{\mH}\geq \e/2\}+\mP\{\|X^n_t-\pi_nh\|_{\mH}\geq r-\e\}\\
&&+\mP\{\|\pi_nh-h\|_{\mH}\geq \e/2\}\\
&<&\eta/2+1-\eta+\eta/2\\
&=&1.
\de
So, the proof is completed.
\end{proof}

In order to present our main result in this section, we need to introduce the following assumption.

\vspace{3mm}

\begin{enumerate}[\bf{Assumption 7.}]
\item
\end{enumerate}
\begin{enumerate}[{\bf (i)}]
\item There exists a Borel measurable mapping $\varrho:\mH\rightarrow\mH$ such that
$$
b(x)=\sigma(x)\varrho(x),
$$
\item
\ce
&&\mE\Big[\exp\Big\{\frac{1}{2}\int_0^T\left\|\varrho(X_s)\right\|_{\mH}^2\dif s
+\int_0^T\int_{\mU_0}\left(\frac{1-\lambda(u)}{\lambda(u)}\right)^2\lambda(u)\nu(\dif u)\dif s\Big\}\Big]\\
&<&\infty.
\de
\end{enumerate}

\vspace{3mm}

Taking
\ce
\Lambda_t:&=&\exp\bigg\{-\int_0^t\<\varrho(X_s),\dif W_s\>_{\tilde{\mH}}-\frac{1}{2}\int_0^t
\left\|\varrho(X_s)\right\|_{\mH}^2\dif s\\
&&\quad\qquad -\int_0^t\int_{\mU_0}\log\lambda(u)N_{\lambda}(\dif s, \dif u)
-\int_0^t\int_{\mU_0}(1-\lambda(u))\nu(\dif u)\dif s\bigg\},
\de
by \cite[Theorem 6 ]{ppks}, we know that $\Lambda_t$ is a exponential martingale under {\bf Assumption 7. (ii)}. Define a new probability measure $\hat{\mP}$ by
$$
\frac{\dif \hat{\mP}}{\dif \mP}=\Lambda_T.
$$
Thus, \cite[Theorem 2.1]{qw2} admits us to obtain that on the new filtered probability space $(\Omega,\mathscr{F},\{\mathscr{F}_t\}_{t\in[0,T]},\hat{\mP})$, $\hat{W}_t:=W_t+\int_0^t\varrho(X_s)\dif s$ is a cylindrical Brownian motion, and the predictable compensator of $N_{\lambda}(\dif t, \dif u)$ is $\dif t\nu(\dif u)$.

Now, we state the main result of this section.

\bt\label{chth}
Suppose that {\bf Assumptions 4-7} are satisfied. Let $v:\mH\rightarrow\mR$ be a $C^2$
scalar function such that $[\nabla v(x)]\in D(A)$ for any $x\in\mH$ and $\|A\nabla v(\cdot)\|_{\mH}$ is
bounded locally, and $\|A\nabla v(\cdot)\|_{\mH}:\mH\to[0,\infty)$ is continuous. Then the Girsanov density $\Lambda_t$ for Eq.\eqref{3}
has the following path-independent property:
\ce
\Lambda_t=\exp\{v(\Gamma)-v(X_t)\}, \quad t\in[0,T],
\de
if and only if
\be
\varrho(x)&=&(\sigma^*\nabla v)(x), \qquad\qquad\qquad\qquad\quad  x\in\mH,\label{e1}\\
\lambda(u)&=&\exp\{v(x+f(x,u))-v(x)\}, \quad (x,u)\in\mH\times\mU_0, \label{e2}
\ee
and $v$ satisfies the following (infinite-dimensional) integro-differential equation,
\be
&&\frac{1}{2}[Tr(\sigma\sigma^*)\nabla^2 v](x)+\frac{1}{2}\|\varrho(x)\|_{\mH}^2
+\langle x,A\nabla v(x)\rangle_{\mH}\no\\
&&+\int_{\mU_0}\Big[e^{v(x+f(x,u))-v(x)}-1-\<f(x,u),\nabla v(x)\>_{\tilde{\mH}}e^{v(x+f(x,u))-v(x)}\Big]\nu(\dif u)=0,
\label{e3}
\ee
where $\sigma^*(x)$ stands for the conjugate of $\sigma(x)$, $\nabla$ and $\nabla^2$ stand for the
first and second {\it Fr\'{e}chet} operators, respectively.
\et

Although the proof of the above theorem is similar to \cite[Theorem 4.3]{qw2}, we prove it here for the readers' convenience.

\begin{proof}
Firstly, let us prove the ``only if" part. By the expression of $\Lambda_t$, it holds that
\ce
\log \Lambda_t&=&-\int_0^t\<\varrho(X_s),\dif W_s\>_{\tilde{\mH}}-\frac{1}{2}\int_0^t
\left\|\varrho(X_s)\right\|_{\mH}^2\dif s-\int_0^t\int_{\mU_0}\log\lambda(u)\tilde{N_\lambda}(\dif s, \dif u)\\
&&-\int_0^t\int_{\mU_0}\(1-\lambda(u)+\lambda(u)\log\lambda(u)\)\nu(\dif u)\dif s.
\de
Besides, \cite[Proposition 1]{qw2} and  {\bf Assumptions 7. (i)} admit us to obtain that
\be
v(\Gamma)-v(X_t)&=&-\int_0^t\langle A X_s,\nabla v(X_s)\rangle_{\mH}\dif s-\int_0^t\langle \varrho(X_s),\sigma^*\nabla v(X_s)\rangle_{\mH}\dif s-\int_0^t\<(\sigma^*\nabla v)(X_s),\dif W_s\>_{\tilde{\mH}}\no\\
&&-\int_{\mU_0}\Big[v(X_{s-}+f(X_{s-},u))-v(X_{s-})-\<f(X_{s-},u),\nabla v(X_{s-})\>_{\tilde{\mH}}\Big]\lambda(u)\nu(\dif u)\dif s\no\\
&&-\int_{\mU_0}\left[v(X_{s-}+f(X_{s-},u))-v(X_{s-})\right]\tilde{N_\lambda}(\dif s, \dif u)\no\\
&&-\frac{1}{2}\int_0^t[Tr(\sigma\sigma^*)\nabla^2 v](X_s)\dif s.
\label{itofor}
\ee
Based on the uniqueness of decomposition for $\log \Lambda_t$ (\cite{dm}), one can have that
\ce
&&\varrho(X_s)=(\sigma^*\nabla v)(X_s), \quad \log\lambda(u)=v(X_{s-}+f(X_{s-},u))-v(X_{s-}),\\
&&-\int_0^t\langle A X_s,\nabla v(X_s)\rangle_{\mH}\dif s-\frac{1}{2}\int_0^t\left\|\varrho(X_s)\right\|_{\mH}^2\dif s+\int_{\mU_0}\<f(X_{s-},u),\nabla v(X_{s-})\>_{\tilde{\mH}}\lambda(u)\nu(\dif u)\dif s\\
&&-\frac{1}{2}\int_0^t[Tr(\sigma\sigma^*)\nabla^2 v](X_s)\dif s+\int_0^t\int_{\mU_0}\(1-\lambda(u)\)\nu(\dif u)\dif s=0.
\de
Thus, by Lemma \ref{supthe}, it holds that (\ref{e1})-(\ref{e3}) are right.

Next, we show ``if" part. Combining (\ref{itofor}) with (\ref{e1})-(\ref{e3}), one can get that
\ce
v(\Gamma)-v(X_t)&=&-\int_0^t\<\varrho(X_s),\dif W_s\>_{\tilde{\mH}}-\frac{1}{2}\int_0^t
\left\|\varrho(X_s)\right\|_{\mH}^2\dif s-\int_0^t\int_{\mU_0}\log\lambda(u)\tilde{N_\lambda}(\dif s, \dif u)\\
&&-\int_0^t\int_{\mU_0}\(1-\lambda(u)+\lambda(u)\log\lambda(u)\)\nu(\dif u)\dif s\\
&=&\log \Lambda_t.
\de
The proof is completed.
\end{proof}

\br
Comparing Theorem \ref{chth} with \cite[Theorem 4.3]{qw2}, one can find that, here $\sigma(x)$ may be degenerate or even could be zero.
\er

The above theorem gives a necessary and sufficient
condition, and hence a characterization of path-independence for the density $\Lambda_t$ of the
Girsanov transformation for a stochastic evolution equation with jumps in terms of an infinite-dimensional
integro-differential equation. Namely, we establish a bridge from Eq.(\ref{3}) to an infinite-dimensional
integro-differential equation.

\bigskip

\textbf{Acknowledgements:}

The authors would like to thank Professor Xicheng Zhang for his fruitful discussions and valuable suggestions.

\end{document}